\input amstex
\magnification=1200
\font\cyr=wncyr10
\font\cyb=wncyb10
\font\cyi=wncyi10
\font\cyre=wncyr8
\font\cyie=wncyi8
\font\cybe=wncyb8
\documentstyle{amsppt}
\NoRunningHeads
\NoBlackBoxes
\define\Cvir{\operatorname{\Bbb Cvir}}
\define\Vir{\operatorname{Vir}}
\define\Ner{\operatorname{Ner}}
\define\Mantle{\operatorname{Mantle}}
\define\Train{\operatorname{Train}}
\define\SI{\operatorname{SI}}
\define\GI{\operatorname{GI}}
\define\BP{\operatorname{BP}}
\define\GM{\operatorname{GM}}
\define\Map{\operatorname{Map}}
\define\FG{\operatorname{FG}}
\define\enl{\operatorname{enl}}
\define\sf{\operatorname{\bold s\bold f}}
\define\cwit{\operatorname{{}^{\circlearrowleft}\frak w\frak i\frak t}}
\define\czwie{\operatorname{{}^{\circlearrowleft}\frak z\frak w\frak i\frak
e}}
\topmatter
\title String field theory and quantum groups. I.
Quantum group structures in geometric quantization of
a self-interacting string field
\endtitle
\author\bf D.V.Juriev\footnote{\ This is an English translation of the
original Russian version, which is located at the end of the article as an
appendix. In the case of any differences between English and Russian
versions caused by a translation the least has the priority as the original
one.\newline}
\endauthor
\date q-alg/9708009\enddate
\affil\eightpoint\it
Research Center for Mathematical Physics and Informatics
"Thalassa Aitheria",\linebreak
ul.Miklukho-Maklaya 20-180, Moscow 117437 Russia\linebreak
E-mail: denis\@juriev.msk.ru
\endaffil
\abstract
The paper is devoted to a description of quantum group structures
in the geometric quantization of a self-interacting string field,
which appear under a transition from a tree-level of the theory to
the account of loop effects in non\-per\-tur\-ba\-tive quantum field
theory of strings.
\endabstract
\endtopmatter
\document
Theory  of strings, which is the most perspective approach to the
unification of Standard Model and its supergeneralizations with quantum
gravity in the unified consistent theory of elementary particles and
their interactions, exists in two forms: as a theory of the first
quantized strings (string quantum mechanics) and as a string field theory
(theory of the second quantized strings) [1,2]. One of the advantages of
the first approach besides the technical circumstances is the clarity
of a geometrical description of the particle interactions, whereas its
essential deficiencies are, first, a difficulty of the consistent account
of the nonperturbative effects, second, the explicit dependence of the
formulation of theory on metric and topology of the background. Both
difficulties obstruct to use the theory as for the problems of quantum
gravity as for the unified quantum description of gravity and Yang-Mills
fields of Standard Model and its supergeneralizations. The second approach
combining string and quantum field ones is free of such disadvantages
(and, more\-over, is of interest for the quantum field theory of vortices in
the quantum fluids), but, its comparative difficulty and awkwardness in
the concerete computations as well as abundance of heterogeneous and
unrelated directly concepts made its applications very inconvenient.
In the author's papers [3] (see also [4]) there was given a unified
formalism for the string field theory based on the geometric quantization,
which allowed to connect other known approaches, for example, the Witten's
polynomial (cubic) string field theory of 1986 [5], the Aref'eva-Volovich
approach, the central place in which belongs to the nonassociative string
algebra [6], the nonpolynomial string theories of B.Zwiebach and other
authors [7]. It was found that despite of the all variability and elegence
of algebraic structure of the first quantized string theory, which includes
the Kac-Moody algebras, the Virasoro algebra $\Cvir$, the Virasoro-Bott
group $\Vir$, the Neretin semigroup $\Ner$, the mantle $\Mantle(\Vir)$ of the
Virasoro-Bott group, the conformal category $\Train(\Vir)$, the train of
this group, and the conformal modular functor among others (see, for example,
[1-3] and references wherein), the algebraic structure of the string field
theory is not less interesting and substantive, and the independence of the
theory itself from metric and topology of the background makes it not only
sometimes simpler than the theory of the first quantized strings (the
multiloop computations or computations on the arbitrary curved background
besides the narrow class of known solutions of the string Einstein
equations are more than labour-consuming), but also permits to account
such effects as the bifurcation of the background or backgrounds with
nontrivial topological and analytical structure of the fractal type (with
infinitely generated fundamental group). Thus, in the second paper from the
unfinished series [3] there was explicated a quantum group structure of the
self-interacting string field in general features (that hypothetically may
connect the known string field approaches with the Vladimirov-Volovich adelic
formalism [8]). This paper is devoted to the more detailed discussion of
the connection between string field theory and the quantum group theory.
The first part contains an analysis of the quantum group structures
in the geometric quantization of the self-interacting string field, whereas
the second part will be devoted to the quantum dynamics of the transition
processes in the string field theory (the so-called ``string cosmological
evolution''). Generally one may said that the quantum group phenomena
naturally appear under a transition from the tree level of the string field
theory to the consistent nonperturbative field description of the loop
effects.

\head\bf 1. Infinite dimensional noncommutative geometry of
a self-interacting string field [3]
\endhead

This paragraph contains an exposition of the main concepts of a string field
the\-ory (bosonic one for simplicity, the supercase needs in slight natural
modifications) in the formalism of geometric quantization applied to a
self-interacting string field (see [3] and also [4]).

The main geometric (geometrodynamic) data of string field theory (as
for a free field as for a self-interacting one) are:\newline
-- The infinite dimensional linear space $\Cal Q$ (or the dual $\Cal Q^*$) of
external degrees of freedom of a string. The coordinates $x^\mu_n$ on $\Cal Q$
are the Taylor coefficients of functions $x^\mu(z)$ that determine the
world-sheet of the string in the complexified target space.\newline
-- The flag manifold of the Virasoro-Bott group $M(\Vir)$ [9] of the
internal degrees of freedom of the string, which is identified, via the
Kirillov construction [10], with the class $S$ of univalent functions $f(z)$;
the natural coordinates on $S$ are coefficients $c_k$ of the Taylor expansion
of the univalent function $f(z)$:
$f(z)=z+c_1z^2+c_2z^3+c_3z^4+\ldots+c_nz^{n+1}+\ldots$.\newline
-- The space $\Cal C$ of the universal deformation of the complex disk [11],
with $M(\Vir)$ as the base and with fibers isomorphic to $D_+$; the
coordinates on $\Cal C$ are $z$, $c_1$, $c_2$, $c_3$, ... $c_n$, ...,
where $c_k$ are the coordinates on the base and $z$ is the coordinate on the
fibers.\newline
-- The space $M(\Vir)\ltimes\Cal Q^*$ of both external and internal degrees of
freedom of the string, or, equivalently, the bundle over $M(\Vir)$
associated with $p:\Cal C\mapsto M(\Vir)$ whose fibers are linear spaces
$\Map(\Cal C/M(\Vir);\Bbb C^n)^*$ dual to the spaces of mappings of fibers of
$p:\Cal C\mapsto M(\Vir)$ into $\Bbb C^n$ (here $\Bbb C^n$ is a local chart
on the background, see [3,4]).\newline
-- The space $\Omega^{\SI}_{\BP}(E_{h,c})$ of {\it semi-infinite
Banks-Peskin differential forms}, which are certain geometric objects on
$M(\Vir)\ltimes\Cal Q^*$ of rather complicated structure [3,4] and
constructed using the prequantization bundle over $M(\Vir)$, here $h$ and $c$
are the prequantization data, in particular $c$ is the central charge
(in general, it differs from the dimension of the background).\newline
-- $Q$ is the natural geometric BRST-operator on
$\Omega^{\SI}_{\BP}(E_{h,c})$; $Q^2=0$ if and only if $c=26$ (with arbitrary
dimension of the background).\newline
-- The space $\Omega_{\BP}^{\SI}(E_{h,c})^*$ of {\it Siegel string fields\/}
with (pseudo)Hermitian metric $(\cdot|\cdot)$.\newline
-- $Q^*$ is the {\it Kato-Ogawa BRST-operator\/} in the space of
Siegel string fields conjugate to $Q$; it defines a new (pseudo)Hermitian
metric $((\cdot|\cdot))=(\cdot|Q^*|\cdot)$ on
$\Omega^{\SI}_{\BP}(E_{h,c})^*$.\newline
-- $\FG_{h,c}(M(\Vir))$ is the {\it Fock-plus-ghost bundle\/} over $M(\Vir)$,
whose sections are just the Banks-Peskin differential forms.\newline
-- The {\it Gauss-Manin string connection\/} $\nabla^{\GM}$ on
$\FG_{h,c}(M(\Vir))$; its covariantly constant sections are the
{\it Bowick-Rajeev vacua}.\newline
-- $D_{\nabla^{GM}}$ is the covariant differential with respect to the
Gauss-Manin string con\-nection. Its nilpotency on the flat background implies
the equality of the central charge to the background's dimension.\newline
-- The space $\Omega^{\SI}_{\BP}(E_{h,c})^*_{\GI}$ og gauge-invariant
Siegel string fields, which is dual to the space of Bowick-rajeev vacua;
this space possess a (pseudo)Hermitian metric $(\cdot|\cdot)_0$, which is
the restriction of the metric $(\cdot|\cdot)$.\newline
-- $Q^*_0$ is the Kato-Ogawa BRST-operator in the space of gauge-invariant
Siegel string fields ($Q^*=D_{\nabla^{\GM}}+Q^*_0$); the (pseudo)Hermitian
metric $((\cdot|\cdot))$ is just the restriction of $((\cdot|\cdot))$ to
$\Omega^{\SI}_{\BP}(E_{h,c})^*_{\GI}$. The existence and nilpotency
of the Kato-Ogawa BRST-operator on the flat background implies the
equality of its dimension to 26.\newline

Equally with the constructed objects with the fixed values of $h$ we shall
consider their direct (discrete or continious) sums over all admissible
$h$ (in particular, $\Omega^{\SI}_{\BP;c}$ as the space of semi-infinite
Banks-Peskin differential forms).

Note that the spaces of Banks-Peskin differential forms, Siegel string
fields, gauge-invariant string fields, Bowick-Rajeev vacua are the
superspaces and objects on them are also superobjects, but we omit
the prefix ``super'' for simplicity. Formulas for the Virasoro algebra
actions in all these spaces in the flat or curved background and for
the BRST-operators are contained in [3,4]. The list of datat above
completely characterizes the free string field theory, the self-interaction
claims to introduce additional algebraic structures. The main spaces of
the theory (the spaces of the Banks-Peskin differential forms and Siegel
string fields) does not depend on the background metric, which determines
as geometrodynamic objects of the string field theory: (pseudo)Hermitian
metrics on the mentioned spaces and BRST-operators, as the gauge
characteristics -- the Gauss-Manin string connection and the related
covariant differential in the Fock-plus-ghost bundle. Note, however, that
the metrics and BRST-operators in the spaces of Banks-Peskin differential
forms and Siegel string fields, which are considered as abstract linear
spaces, in\-dis\-tin\-guishable for various backgrounds (the independence
of the second quantized string field theory from the background), the metric
on the background is restored only under the consideration of the spaces
as spaces of geometric objects on the space of both external and internal
degrees of freedom of a string. It allows to give a traditional for the
string field theory interpretation of metrics on the background and Yang-Mills
fields as ``low-temperature limits'' of fields of closed and open strings
(though the alternatives are possible, for example, to consider the
components of Higgs fields for the Yang-Mills fields of Standard Model in
the same sector as a gravitational one, i.e. in the closed string sector,
in view of the presence of the subsidiary nonmetric degrees of freedom of
cohomology of the Virasoro algebra with coefficients in string fields [4,3]).
The flatness conditions for the Gauss-Manin connection (or, equivalently,
for the nilpotency of the covariant differential) are just the {\it string
Einstein equations\/} (which transform into ordinary ones in the
``low-temperature limit''). Therefore, the string Einstein equations may be
also defined as conditions of the existence and the nilpotency of the
Kato-Ogawa BRST-operator in the space of gauge-invariant Siegel string fields.

If the background does not obey the string Einstein equations then the
Bowick-Rajeev vacua do not exist, there are some ways to construct
string field theory in such situation, for example, to use the
Bowick-Rajeev instantons [3]. The questions of (in)dependence of the
resulted theories from the background as on the connected components of
the space of solutions of the string Einstein equations as in whole are
discussed in the second paper from the series [3] (see also refs wherein).

To formulate the self-interacting string field theory it is convenient to
use the formalism of the local conformal field algebras, which is based on
the ideas of the noncommutative geometry and which is exposed in details
in the original papers [12] and a review [13] (see also [14]).

Let us consider as in [3] the space $\Omega^{\cdot}_{\BP;\enl}=
\Omega^{\cdot}(\tilde\Bbb C^*,\Omega^{\SI}_{\BP;c})$ of enlarged
Banks-Peskin differential forms (here $\tilde\Bbb C^*$ is the universal
covering of the complex plane punctured at zero). Let us define also
the space of the enlarged Siegel string fields $\Omega^{\cdot}_{\sf;\enl}
=\Omega^{\cdot}(\tilde\Bbb C^*,(\Omega^{\SI}_{\BP;c})^*)$. Formulas for
the action of the Virasoro algebra in the spaces $\Omega^{\cdot}_{\BP;\enl}$
and $\Omega^{\cdot}_{\sf;\enl}$ are contained in [3]. Let us construct also
the enlarged BRST-operators $Q_{\enl}$ and $Q^*_{\enl}$ as exterior
differentials in the spaces $\Omega^{\cdot}_{\BP;\enl}$ and
$\Omega^{\cdot}_{\sf;\enl}$ from the BRST-operators $Q$ and $Q^*$.

\proclaim{Theorem 1 [3]} The space $\Omega^{\cdot}_{\sf;\enl}$ admits the
structure of a BRST-invariant local conformal field algebra that is
covariant with respect to the Gauss-Manin connection $\nabla_{\GM}$.
\endproclaim

Thus, the space $\Omega^{\cdot}_{\sf;\enl}$ may be regarded as a
noncommutative de Rham complex (cf.[15,16]) with respect to the
enlarged BRST-operators. This complex is called the {\it enlarged string
field algebra}. Relations of the enlarged string field algebra to
the {\it Aref'eva-Volovich nonassociative string field algebra}, which
is realized in the space of Siegel string fields and which is a certain
reduction of the associative enlarged string field algebra, are described
in [3].

The elements of the enlarged string field algebra $\Omega^{\cdot}_{\sf;\enl}$
form a Lie algebra under the commutator. This Lie algebra admits a central
extension by the imaginary part of the (pseudo)Hermitian metric
$((\cdot|\cdot))$. Consider the connection forms on $\tilde\Bbb C^*$ with
values in $(\Omega^{\SI}_{\BP;c})^*$, i.e. gauge fields on $\tilde\Bbb C^*$
valued in Siegel string fields; elements of $\Omega^0_{\sf;\enl}$ realize
the infinitesimal gauge transformations of these fields. These gauge
transformations are closed (so we are in the situation of the Witten's
string field theory of 1986 [5]), the corresponding Lie algebra is called
the {\it Witten string Lie algebra\/} and is denoted by $\cwit$ (the circled
arrow $\circlearrowleft$ is the code for ``string''). The space of
$\nabla_{\GM}$--covariant elements of the Lie algebra $\cwit$ is denoted by
$\cwit_{\nabla^{\GM}}$ and is also called the witten string Lie algebra.
The Witten string Lie algebra $\cwit$ is just the central extension of
the commutator algebra of the zero component of the enlarged string field
algebra, which is described above.

The are canonical (Lie-Berezin) Poisson brackets on the space $\cwit^*$ (or
$\cwit^*_{\GM}$) dual to the Witten string Lie algebra $\cwit$ (or
$\cwit_{\GM}$), which can be quantized as such.

The Lie-Berezin brackets in the coadjoint representation of the Witten
string Lie algebra may be reduced to the nonpolynomial brackets in the
space of all functionals on the Banks-Peskin differential forms (or
Bowick-Rajeev vacua), a procedure of the Hamiltonian reduction is
described in [3] and follows the general scheme of reduction of Lie-Berezin
brackets (see, for example, [17]). These nonpolynomial brackets generate
a Lie quasi(pseudo)algebra (quasialgebra in the terminology of [18], see
also [19], and pseudoalgebra in the terminology of [17]) of infinitesimal
nonpolynomial gauge transformations. The nonpolynomial transformations
on the space of Bowick-Rajeev vacua were considered in [7]; they generate
a Lie quasi(pseudo)algebra, which is denoted by $\czwie_{\nabla^{\GM}}$
and is called the Zwiebach string Lie quasi(pseudo)algebra; the related
quasi(pseudo)algebra on the space of Banks-Peskin differential forms
is denoted by $\czwie$ and has the same name. Nonpolynomial brackets
are realized in functionals on the space $\czwie^*$ (or
$\czwie^*_{\nabla^{\GM}}$) dual to the Lie quasi(pseudo)algebra $\czwie$
(or $\czwie_{\nabla^{\GM}}$). The Zwiebach string Lie quasi(pseudo)algebra
may be obtained from the Aref'eva-Volovich nonassociative string field
algebra as its ``commutator'' algebra. More precisely, the higher
ope\-rations of the Sabinin-Mikheev multialgebra [20] constructed from the
Zwiebach string Lie quasi(pseudo)algebra are just the higher commutators
in the Aref'eva-Volovich nonassociative string field algebra.

Thus, the nonpolynomial string field theory [7] in the space of
Banks-Peskin differential forms (or Bowick-Rajeev vacua) can be obtained
from the cubic Witten-type string field theory [5] in the enlarged
space by use of the Hamiltonian reduction. Moreover, the approach of the
papers [7] on the nonpolynomial field theory appears to be equivalent
to the approach of I.Ya.Aref'eva and I.V.Volovich [6] based on the
nonassociative string field algebra.

\head\bf 2. Quantum group structure of a self-interacting string field
\endhead

Note that there are two ways to quantize the self-interacting string field.
First, one may quantize the nonpolynomial Poisson brackets themselves,
for instance in the formalism of asymptotic quantization [17]. Second,
quantum theory may be obtained by a quantum reduction of the quantized
Lie-Berezin brackets on the space $\cwit_{\nabla^{\GM}}$; in this case
the algebra of quantum observables is identified with certain quantum
reduction of the universal envelopping algebra $\Cal U(\cwit_{\nabla^{\GM}})$
of the Witten string Lie algebra. Both variants are of a undoubtful
interest from the mathematical point of view and look rather natural.
However, below we shall try to unravel a very important nuance significantly
changing the initial ``na\"\i ve'' point of view.

We stress that the objects constructed above describe a self-interacting
string field theory only on the tree level, i.e. define the so-called
``classical string field theory'' in the terminology of [21]. To describe
the string field theory completely and consistently in the nonperturbative
mode it is necessary to use the following crucial result.

\proclaim{Theorem 2 [3]} The Witten string Lie algebra $\cwit_{\nabla^{\GM}}$)
(or $\cwit$) admits the structure of a Lie bialgebra.
\endproclaim
\rm

Proof of this theorem is based on the auxiliary statement that the
enlarged string field algebra is a crossing-algebra [3].

Thus, there is explicated a quantum group structure of a self-interacting
string field theory on the quasiclassical level (cf.[22,17]). So on the
quantum level [23] the algebra of observables is described by the quantum
universal envelopping algebra $\Cal U_q(\cwit_{\nabla^{\GM}})$ (or
$\Cal U_q(\cwit)$), or, precisely, by some its quantum reduction,
however, the explicit construction of such infinite dimensional Hopf algebra
is not known. Before to pass to a description of the quantum algebra of
observables for the self-inteacting string field theory in the concrete
cases let us examine the process of reduction on the quasiclassical level,
otherwords, determine in what the Witten string Lie bialgebra is
transformed under the reduction of the Lie-Berezin brackets constructed
from the Witten string Lie algebra to the nonpolynomial Poisson brackets
(remind, that the Witten string Lie algebra itself transforms into
the Zwiebach string Lie quasi(pseudo)algebra).

\proclaim{Theorem 3} The Zwiebach string Lie quasi(pseudo)algebra
$\czwie_{\nabla^{\GM}}$ (or $\czwie$) admits the structure of
the cojacobian quasibialgebra.
\endproclaim
\rm

Cojacobian quasibialgebras [24] form a class of Lie protobialgebras
dual to the jacobian quasibialgebras (Lie quasibialgebras in the
terminology of V.G.Drinfeld [25]).

To prove the theorem one should apply a reduction to the double of
the witten string Lie bialgebra with a translation invariant bracket.

Remind that the cojacobian quasibialgebras are the infinitesimal objects
for the Poisson quasigroups [24] (whereas the jacobian quasibialgebras --
for quasiPoisson Lie groups [25]). Thus, the Zwiebach string Lie
quasibialgebra realizes a qua\-si\-clas\-si\-cal version of the nonlinear
geometric algebra [26] on the infinitesimal level, whose purely quantum
version is not known nowadays. A relation between the structure of
cojacobian quasibialgebra and such object of the nonlinear geometric
algebra as Sabinin-Mikheev multialgebra was explicated in [24].
A specific character of the infinite dimensional situation is in the fact
that perhaps none global quasigroup corresponds to the mentioned
infinitesimal objects.

Let us pass now from the quasiclassics to the explicit construction of
the (en\-larged) quantum algebra of observables
$\Cal U_q(\cwit_{\nabla^{\GM}})$
in the concerete cases (and without acoount of ghosts, i.e. in the Fock
sector, for simplicity). Consider a flat compact background isomorphic to
the quotient of the euclidean space by the root lattice of the semisimple
Lie algebra $\frak g$. In this case the enlarged string field algebra in
the space of enlarged gauge-invariant Siegel string fields
$\Omega^0_{\sf;\enl}$ (more precisely, in its Fock sector
$\Omega^0_{\sf;\enl;F}$) is a local conformal field algebra received by
a renormalization of the pointwise product of operator fields [3,13,14] (see
also [27]) from the vertex operator algebra constructed by this lattice [28].
Hence, the linear space $\Omega^0_{\sf;\enl;F}$ may be identified with the
space $\Cal U(\hat\frak g_+))[[t]]$ of semi-infinite formal power series
with coefficients in the universal envelopping algebra of the positive
(nonnegative) component of the Kac-Moody algebra $\frak g_+$. Note that
the vertex operator algebra is generated by its currents (primary fields
of spin 1), the components of which from the Kac-Moody algebra $\hat\frak g$,
therefore the enlarged string field algebra is a quotient of the universal
envelopping algebra $\Cal U(\hat\frak g)$ of the Lie algebra $\hat\frak g$ by
its ideal $\Cal J$. Therefore, the Witten string Lie algebra
$\cwit_{\nabla^{\GM};F}$ (the symbol $F$ means the Fock sector) is a
quotient of the commutator algebra $\Cal U_{[\cdot,\cdot]}(\hat\frak g)$
by the ideal $\Cal J_{[\cdot,\cdot]}$. The quantum version of the Witten
string Lie algebra can be obtained in the following way: consider the
quantum universal envelopping algebra $\Cal U_q(\hat\frak g)$ supplied by
the $q$-commutator; in view of an existence of the $q$-vertex construction
for this algebra, the ideal $\Cal J$ can be deformed to an ideal $\Cal J_q$
of the algebra $\Cal U_q(\hat\frak g)$, which is closed under the
$q$-commutator; the relations between the elements of
$\Cal U_q(\hat\frak g)/\Cal J_q$ defined by $q$-commutator are just the
defining ones in the quantum Witten string algebra
$\Cal U_q(\cwit_{\nabla^{\GM};F})$.

Note taht in this important for the string theory class of examples
on the intermediate step of the construction of the quantum Witten string
algebra realizing a quantum version of gauge symmetries of string fields
an object (enlarged string field algebra), which is natural to be considered
from point of view of the 2-loop formalism in theories of strings and
integrable systems [29], appeared. Moreover, it will be too interesting to
examine how the modular invariance of the first quantized string explicates
itself in the quantum group formalism of the nonperturbative string field
theory.

Thus, various quantum group structures of the string field theory as
on the quasiclassical as on the quantum levels were explored in the paper.
General state\-ments were formulated, examples were examined, new aspects
of the theory of quantum groups were demonstrated (i.e. relations to the
nonlinear geometric algebra, namely, the quasiclassical versions of
so-called ``quantum quasigroups'' and ``quantum loops'') in the context of
the second quantization of a string under a transition from the tree level
to the consistent account of the loop effects in nonperturbative string
field theory.

\Refs
\roster
\item"[1]" {\it Green M., Schwarz J., Witten E.}, Superstring theory.
Cambridge Univ.Press, Cambridge, 1987.
\item"[2]" {\it Morozov A.Yu., Perelomov A.M.}, Current Probl.Math.,
Fundam.Directions. M., VINITI, 1989;\newline
{\it Morozov A.Yu.}, Elem.Part.Atom.Nuclei 23(1) (1992) 174-238;\newline
{\it Morozov A.Yu.}, Soviet Phys.Uspekhi 35 (1992) 671-714.
\item"[3]" {\it Juriev D.}, Alg.Groups Geom. 11 (1994) 145-179 [e-version:
hep-th/9403068]; Russian J.Math.Phys. 4 (1996) 287-314; J.Geom.Phys. 16 (1995)
275-300.
\item"[4]" {\it Juriev D.}, Lett.Math.Phys. 22 (1991) 1-6, 11-14;\newline
{\it Juriev D.}, Lett.Math.Phys. 19 (1990) 355-356; 19 (1990) 59-64.
\item"[5]" {\it Witten E.}, Nucl.Phys.B268 (1986) 253-294.
\item"[6]" {\it Aref'eva I.Ya., Volovich I.V.}, Phys.Lett.B182 (1986)
159-163, 312-316; 189 (1987) 488.
\item"[7]" {\it Saadi M., Zwiebach B.}, Ann.Phys.(NY) 192 (1989)
213-227;\newline
{\it Kugo T., Kunimoto H., Suehiro K.}, Phys.Lett.B226 (1989) 48-54;\newline
{\it Sonoda H., Zwiebach B.}, Nucl.Phys.B331 (1990) 592-628;\newline
{\it Kugo T., Suehiro K.}, Nucl.Phys.B337 (1990) 434-466;\newline
{\it Zwiebach B.}, Commun.Math.Phys. 136 (1991) 83-118.
\item"[8]" {\it Volovich I.V.}, Class.Quant.Grav. 1987. V.4. P.83;\newline
{\it Vladimirov V.S.}, Lett.Math.Phys. 1993. V.27. P.123.
\item"[9]" {\it Segal G.}, Commun.Math.Phys. 80 (1981) 301-342;\newline
{\it Kirillov A.A.}, Lect.Notes Math. 970 (1984) 49-67;\newline
{\it Kirillov A.A., Juriev D.V.}, Funct.Anal.Appl. 20 (1986) 322-324;
21 (1987) 284-293;\newline
{\it Witten E.}, Commun.Math.Phys. 114 (1988) 1-53;\newline
{\it Kirillov A.A.}, in ``Infinite-dimensional Lie algebras and quantum field
theory''. World Scientific, Teaneck, NJ, 1988, pp.73-77;\newline
{\it Kirillov A.A., Juriev D.V.}, J.Geom.Phys. 5 (1988) 351-364;\newline
{\it Kirillov A.A.}, in ``Operator algebras, unitary representations,
envelopping algebras and invariant theory''. Birkhauser, Boston, 1992,
pp.73-83;\newline
{\it Juriev D.}, Adv.Soviet Math. 2 (1991) 233-247;\newline
{\it Kirillov A.A.}, Contemp.Math. 145 (1993) 33-63.
\item"[10]" {\it Kirillov A.A.}, Funct.Anal.Appl. 21(2) (1987) 122-125.
\item"[11]" {\it Juriev D.V.}, St.-Petersburg Math.J. 2 (1991) 401-417;
{\it Juriev D.}, Russian J.Math.Phys. 2 (1994) 111-121.
\item"[12]" {\it Juriev D.}, Commun.Math.Phys. 138 (1991) 569-581,
1992. V.146. P.427; J.Funct.Anal. 1991. V.101. P.1-9.
\item"[13]" {\it Juriev D.V.}, Russian Math.Surveys 46(4) (1991) 135-163.
\item"[14]" {\it Juriev D.V.}, Theor.Math.Phys. 101 (1994) 1387-1403.
\item"[15]" {\it Connes A.}, Introduction \`a la g\'eom\'etrie non commutative.
InterEditions, Paris, 1990.
\item"[16]" {\it Manin Yu.I.}, Topics in non-commutative geometry. Princeton
Univ.Press, Princeton, NJ, 1991.
\item"[17]" {\it Karasev M.V., Maslov V.P.}, Nonlinear Poisson brackets.
geometry and quantization. Amer.Math.Soc., Providence, RI, 1993.
\item"[18]" {\it Batalin I.}, J.Math.Phys. 22 (1981) 1837-1850.
\item"[19]" {\it Mikheev P.O.}, in ``Some applications of differential
geometry''. Moscow, 1985, pp.85-93 [VINITI: 4531-85Dep.];\newline
{\it Mikheev P.O.}, Trans.Inst.Phys. Estonian Acad.Sci. 66 (1990) 54-66.
\item"[20]" {\it Sabinin L.V., Mikheev P.O.}, Soviet Math. 36 (1988) 545-548.
\item"[21]" {\it Zwiebach B.}, Nucl.Phys.B390 (1993) 33-152.
\item"[22]" {\it Drinfeld V.G.}, DAN SSSR 268(2) (1983) 285-287;\newline
{\it Semenov-Tian-Shanskii M.A.}, Funct.Anal.Appl. 17 (1983) 259.
\item"[23]" {\it Drinfeld V.G.}, Proc.Intern.Congr.Math. (1986), Berkeley,
California, vol.1, pp.789-820;\newline
{\it Reshetihin N.Yu., Takhtadzhan L.A., Faddeev L.D.},
 St.-Petersburg Math.J. 1 (1990) 193-225;\newline
{\it Isaev A.P.}, Elem.Part.Atom.Nuclei 26 (1995) 1204-1263;\newline
{\it Zhelobenko D.P.}, Representations of reductive Lie algebras. Moscow,
Nauka, 1994.
\item"[24]" {\it Bangoura M.}, C.R.Acad.Sci.Paris I 319 (1994) 975-978;\newline
{\it Bangoura M.}, Th\`ese l'Univ.Sci.Tech.Lille (1995) no.1387.
\item"[25]" {\it Drinfeld V.G.}, St.-Petersburg Math.J. 1 (1990)
1419-1457;\newline
{\it Kosmann-Schwarzbach Y.}, C.R.Acad.Sci.Paris I 312 (1991) 391-394;\newline
{\it Kosmann-Schwarzbach Y.}, in ``Mathematical aspects of field theory''.
Amer.Math.Soc., Providence, RI, 1992, pp.459-489.
\item"[26]" {\it Sabinin L.V.}, in ``Webs and quasigroups''. Kalinin
[Tverp1], 1988, pp.32-37;\newline
{\it Sabinin L.V.}, Methods of nonassociative algebra in differential
geometry. Supplement to Russian transl.of Kobayashi S., Nomidzu K.,
Foundations of differential geometry. V.1. Moscow, Nauka, 1982;\newline
{\it Sabinin L.V.}, Trans.Inst.Math.Siberian Branch USSR Acad.Sci. 14 (1989)
208-221;\newline
{\it Mikheev P.O., Sabinin L.V.}, in ``Quasigroups and loops. Theory and
applications''. Ber\-lin, Heldermann Verlag, 1990, pp.357-430.
\item"[27]" {\it Juriev D.V.}, Theor.Math.Phys. 93 (1992) 1101-1105.
\item"[28]" {\it Frenkel I., Lepowsky J., Meurman A.}, Vertex operator
algebras and the Monster. Acad. Press, New York, 1988;\newline
{\it Frenkel I.B., Huang Y.-Z., Lepowsky J.}, Memoires AMS. 1993.
V.494;\newline
{\it Dong C., Lepowsky J.}, Generalized vertex algebras and relative vertex
operators. Birk\-h\"au\-ser, 1993;\newline
{\it Kac V.G.}, Vertex operator algebras for beginners. Birkh\"auser, 1996.
\item"[29]" {\it Morozov A.Yu.}, Soviet Phys.Uspekhi 37(1) (1994) 1-109.
\endroster
\endRefs

\newpage
\centerline{\bf APPENDIX: THE ORIGINAL RUSSIAN VERSION OF ARTICLE}
\ \newline
\centerline{\cyb STRUNNAYa TEORIYa POLYa I KVANTOVYE GRUPPY. {\bf I}.}
\centerline{\cyb KVANTOVOGRUPPOVYE STRUKTURY V GEOMETRIChESKOM}
\centerline{\cyb KVANTOVANII SAMODEI0STVUYuShChEGO STRUNNOGO POLYa}
\ \newline
\centerline{\cyb D.V.Yurp1ev}
\ \newline\eightpoint
\centerline{\cyre Tsentr matematicheskoe0 fiziki i informatiki ``Talassa
E1teriya'',}
\centerline{\cyre ul.Miklukho-Maklaya 20-180, Moskva 117437 Rossiya.}
\centerline{E-mail: denis\@juriev.msk.ru}
\ \tenpoint\newline
\centerline{q-alg/9708009}
\ \newline

\ \newline

\eightpoint\cybe REZYuME. \cyre
Dannaya rabota posvyashchena opisaniyu kvantovogruppovyh struktur v
geo\-met\-ri\-ches\-kom kvantovanii samodee0stvuyushchego strunnogo polya,
voznikayushchih pri pe\-re\-ho\-de s drevesnogo urovnya teorii k uchetu
petlevyh e1ffektov v neperturbativnoe0 kvan\-to\-vo\-po\-le\-voe0 teorii
strun.

\tenpoint\ \newline

\ \newline

\cyr
Teoriya (super)strun, predstavlyayushchaya soboe0 naibolee
sovremennye0 i perspektivnye0 podhod k obp2edineniyu Standartnoe0
Modeli ili ee su\-per\-obob\-shche\-nie0 s kvantovoe0 gravitatsiee0 v
edinuyu posledovatelp1nuyu teoriyu e1lementarnyh chastits i ih
vzaimodee0stvie0, cushchestvuet v dvuh vidah: kak teoriya
pervichnokvantovannyh strun (strunnaya kvantovaya mehanika) i
kak strunnaya teoriya polya (teoriya vtorichnokvantovannyh strun)
[1,2]. K chis\-lu dostoinstv pervogo podhoda naryadu s tehnicheskimi
obstoyatelp1stvami otnosit\-sya prozrachnostp1 geometricheskogo
opisaniya protsessov vza\-i\-mo\-dee0\-s\-t\-viya chastits, v to vremya kak
ee sushchestvennymi nedostatkami yavlyayut\-sya, vo-pervyh,
trudnostp1 posledovatelp1nogo ucheta neperturbativnyh e1ffektov,
vo-vtoryh, yavnaya zavisimostp1 formulirovki teorii ot metriki i
to\-po\-lo\-gii be1kgraunda. Obe trudnosti prepyat\-stvuyut ispolp1zovaniyu
e1toe0 teo\-rii kak dlya zadach kvantovoe0 gravitatsii, tak i dlya
edinogo kvantovogo opi\-sa\-niya gravitatsii i polee0 Yanga-Millsa
Standartnoe0 Modeli i ee su\-per\-obob\-shche\-nie0. Vtoroe0 podhod v teorii
strun, strunnaya teoriya polya, kom\-bi\-ni\-ru\-yu\-shchie0 strunnye0 i
kvantovopolevoe0 podhody, svoboden ot uka\-zan\-nyh nedostatkov (i, krome
togo, predstavlyaet interes dlya kvan\-to\-vo\-po\-le\-voe0 teorii vihree0 v
kvantovyh zhidkostyah), odnako, ego sravnitelp1naya slozhnostp1 i
gromozdkostp1 pri konkretnyh vychisleniyah, a takzhe obilie
raznorodnyh i vneshne ne svyazannyh mezhdu soboe0 kontseptsie0
delali ego primenenie do nedavnego vremeni zatrudnitelp1nym. V
rabotah avtora [3] (sm.takzhe [4]) byl dan edinye0 formalizm
strunnoe0 teorii polya, os\-no\-vy\-va\-yu\-shchie0\-sya na geometricheskom
kvantovanii i pozvolivshie0 svyazatp1 mezhdu soboe0 drugie izvestnye
podhody, naprimer, baziruyushchuyusya na ne\-kom\-mu\-ta\-tiv\-noe0 geometrii
polinomialp1nuyu (kubichnuyu) strunnuyu teoriyu polya Vittena 1986
goda [5], podhod Arefp1evoe0--Volovicha, tsentralp1noe mesto v
kotorom zanimaet neassotsiativnaya strunnaya algebra [6],
ne\-po\-li\-no\-mi\-alp1\-nye strunnye teorii Tsvibaha i drugih avtorov [7].
Pri e1tom okazalosp1, chto pri vsee0 mnogogrannosti i izyashchestve
algebraicheskoe0 struktury teo\-rii pervichnokvantovannyh strun,
vklyuchayushchee0 v sebya algebry Katsa--Mu\-di, algebru Virasoro
$\Cvir$, gruppu Virasoro--Botta $\Vir$, polugruppu Ne\-re\-ti\-na
$\Ner$, mantiyu $\Mantle(\Vir)$ gruppy Virasoro--Botta, konformnuyu
ka\-te\-go\-riyu $\Train(\Vir)$, shlee0f e1toe0 gruppy, i konformnye0
modulyarnye0 funk\-tor sredi prochih (sm.naprimer, [1-3] i ssylki
v nih), algebraicheskaya struktura strunnoe0 teorii polya ne menee
interesna i soderzhatelp1na, a ne\-za\-vi\-si\-mostp1 samoe0 teorii ot
metriki i topologii be1kgraunda delaet ee ne tolp1ko inogda bolee
prostoe0, chem teoriya pervichnokvantovannyh strun (mno\-go\-pet\-le\-vye
vychisleniya ili vychisleniya na obshchem iskrivlennom be1k\-gra\-un\-de
vne sravnitelp1no uzkogo klassa izvestnyh reshenie0 strunnyh
ura\-v\-ne\-nie0 E1e0nshtee0na v kotoroe0 bolee chem trudoemki), no i
pozvolyayushchee0 uchi\-ty\-vatp1 e1ffekty tipa perestroek be1kgraunda
ili be1kgraundy s ne\-tri\-vi\-alp1\-noe0 topologicheskoe0 i analiticheskoe0
strukturoe0 fraktalp1nogo tipa (beskonechnoporozhdennoe0
fundamentalp1noe0 gruppoe0). Tak vo vtoroe0 rabote iz neokonchennogo
tsikla [3] byla vyyavlena v obshchih chertah kvan\-to\-vo\-grup\-po\-vaya
struktura samodee0stvuyushchego strunnogo polya (chto,
pred\-po\-lo\-zhi\-telp1\-no,
mozhet svyazyvatp1 imeyushchiesya strunnopolevye podhody s adelp1nym
formalizmom Vladimirova--Volovicha [8] i s ideyami kvan\-to\-va\-niya
prostranstva-vremeni v tselom, sm.naprimer, [9:\S\S33,34]). Bolee
detalp1nomu obsuzhdeniyu svyazi strunnoe0 teorii polya s teoriee0
kvantovyh grupp i posvyashchena dannaya rabota. V pervoe0 chasti
obsuzhdayut\-sya kvan\-to\-vo\-grup\-po\-vye struktury v geometricheskom
kvantovanii samodee0stvuyushchego strunnogo polya, v to vremya
kak vtoraya chastp1 budet posvyashchena kvantovoe0 dinamike
perehodnyh protsessov v strunnoe0 teorii polya (t.n. ``strunnoe0
kosmologicheskoe0 e1volyutsii''). V tselom, mozhno skazatp1, chto
kvan\-to\-vo\-grup\-po\-vye yavleniya voznikayut estestvennym obrazom pri
perehode s dre\-ves\-no\-go urovnya v strunnoe0 teorii polya k
posledovatelp1nomu ne\-per\-tur\-ba\-tiv\-no\-mu polevomu opisaniyu petlevyh
e1ffektov.

\head\cyb 1. Beskonechnomernaya nekommutativnaya geometriya
samodee0stvuyushchego strunnogo polya [3]
\endhead

\cyr
V dannom paragrafe kratko izlagayut\-sya osnovnye ponyatiya strunnoe0
teorii polya (dlya prostoty bozonnoe0, supersluchae0 trebuet
ne\-zna\-chi\-telp1\-nyh estestvennyh modifikatsie0) v formalizme
geometricheskogo kvan\-to\-va\-niya primenitelp1no k
samodee0stvuyushchemu strunnomu polyu (sm.[3], a takzhe [4]).

Osnovnymi geometricheskimi (geometrodinamicheskimi) dannymi strun\-noe0
teorii polya (kak dlya svobodnogo polya, tak i dlya samodee0stvuyushchego)
yavlyayut\-sya:\newline
-- Beskonechnomernoe linee0noe prostranstvo $\Cal Q$ (ili emu
dvoe0stvennoe $\Cal Q^*$) vneshnih stepenee0 svobody struny.
Koordinaty $x^\mu_n$ na $\Cal Q$ sutp1 tee0\-lo\-rov\-s\-kie
koe1ffitsienty funktsie0 $x^\mu(z)$, opredelyayushchih mirovoe0 list
struny v kompleksifitsirovannom prostranstve.\newline
-- Mnogoobrazie flagov gruppy Virasoro-Botta $M(\Vir)$ [10] vnutrennih
stepenee0 svobody struny, otozhdestvlyaemoe posredstvom konstruktsii
Ki\-ril\-lo\-va [11] s klassom $S$ odnolistnyh funktsie0 $f(z)$; estestvennye
ko\-or\-di\-na\-ty na $S$ -- koe1ffitsienty $c_k$ tee0lorovskogo razlozheniya
odnolistnoe0 funktsii $f(z)$:
$f(z)=z+c_1z^2+c_2z^3+c_3z^4+\ldots+c_nz^{n+1}+\ldots$.\newline
-- Prostranstvo $\Cal C$ universalp1noe0 deformatsii kompleksnogo diska
[12] s $M(\Vir)$ v kachestve bazy i sloyami, izomorfnymi edinichnomu
kompleksnomu disku $D_+$; koordinaty na $\Cal C$ sutp1 $z$, $c_1$,
$c_2$, $c_3$, ... $c_n$, ..., gde $c_k$ -- koordinaty na baze i $z$ --
koordinata v sloyah (prostranstvo universalp1noe0 de\-for\-ma\-tsii kak
rassloenie dopuskaet estestvennuyu trivializatsiyu).\newline
-- Prostranstvo $M(\Vir)\ltimes\Cal Q^*$ kak vneshnih, tak i
vnutrennih stepenee0 svo\-bo\-dy struny, ili, chto e1kvivalentno,
rassloenie nad $M(\Vir)$, as\-so\-tsi\-i\-ro\-van\-noe s
$p:\Cal C\mapsto M(\Vir)$, chp1i sloi -- linee0nye prostranstva
$\Map(\Cal C/M(\Vir);\Bbb C^n)^*$, dvoe0stvennye k prostranstvam
otobrazhenie0 sloev rassloeniya $p:\Cal C\mapsto M(\Vir)$ v $\Bbb C^n$
(zdesp1 $\Bbb C^n$ -- lokalp1naya karta na be1kgraunde, sm.[3,4]).\newline
-- Prostranstvo $\Omega^{\SI}_{\BP}(E_{h,c})$ {\cyi polubeskonechnyh
differentsialp1nyh form Be1nksa--Peskina}, e1ti formy sutp1 nekotorye
differentsialp1nye obp2ekty na $M(\Vir)\ltimes\Cal Q^*$ dovolp1no
slozhnoe0 struktury [3,4], svyazannye s rassloeniem pred\-kvan\-to\-va\-niya
na $M(\Vir)$, gde $h$ i $c$ -- dannye predkvantovaniya, v
chastnosti $c$ -- tsent\-ralp1\-nye0 zaryad (voobshche govorya,
otlichnye0 ot razmernosti be1k\-gra\-un\-da).\newline
-- $Q$ -- estestvennye0 geometricheskie0 BRST-operator na
$\Omega^{\SI}_{\BP}(E_{h,c})$; $Q^2=0$, esli i tolp1ko esli $c=26$
(pri e1tom razmernostp1 be1kgraunda pro\-iz\-volp1\-na).\newline
-- Prostranstvo $\Omega_{\BP}^{\SI}(E_{h,c})^*$ {\cyi strunnyh polee0
Zigelya\/} s (psevdo)e1rmitovoe0 metrikoe0 $(\cdot|\cdot)$.\newline
-- $Q^*$ -- {\cyi BRST-operator Kato--Ogavy\/} v prostranstve
strunnyh polee0 Zi\-ge\-lya, sopryazhennye0 $Q$; on opredelyaet
novuyu (psevdo)e1rmitovu metriku $((\cdot|\cdot))=(\cdot|Q^*|\cdot)$
na $\Omega^{\SI}_{\BP}(E_{h,c})^*$.\newline
-- $\FG_{h,c}(M(\Vir))$ -- {\cyi fokovsko-gostovoe rassloenie\/}
nad $M(\Vir)$, chp1i secheniya -- polubeskonechnye differentsialp1nye
formy Be1nksa--Peskina.\newline
-- {\cyi Strunnaya svyaznostp1 Gaussa--Manina\/} $\nabla^{\GM}$ v
rassloenii $\FG_{h,c}(M(\Vir))$, kovariantno postoyannye secheniya
fokovsko-gostovogo rassloeniya sutp1 {\cyi va\-ku\-u\-my Bovika--Radzhiva}.
\newline
-- Kovariantnye0 differentsial $D_{\nabla^{GM}}$ -- kovariantnye0
differentsial po otnosheniyu k svyaznosti Gaussa--Manina. Na ploskom
be1kgraunde ego nilp1\-po\-tent\-nostp1 vlechet ravenstvo razmernosti
be1kgraunda i tsentralp1nogo zaryada.\newline
-- Prostranstvo $\Omega^{\SI}_{\BP}(E_{h,c})^*_{\GI}$
kalibrovochno-invariantnyh strunnyh polee0 Zigelya,
dvoe0stvennoe k prostranstvu vakuumov Bovika--Radzhiva; e1to
pro\-st\-ran\-st\-vo nadelyaet\-sya (psevdo)e1rmitovoe0 metrikoe0
$(\cdot|\cdot)_0$, reduktsiee0 met\-ri\-ki
$(\cdot|\cdot)$.\newline
-- $Q^*_0$ -- BRST-operator Kato-Ogavy v prostranstve
ka\-lib\-ro\-voch\-no-in\-va\-ri\-ant\-nyh strunnyh polee0 Zigelya
($Q^*=D_{\nabla^{\GM}}+Q^*_0$); (psevdo)e1rmitova metrika $((\cdot|\cdot))$ --
reduktsiya $((\cdot|\cdot))$ na $\Omega^{\SI}_{\BP}(E_{h,c})^*_{\GI}$.
Sushchestvovanie i nilp1\-po\-tent\-nostp1 BRST-operatora Kato-Ogavy na
ploskom be1kgraunde vlechet ravenstvo ego razmernosti chislu 26.\newline

Naryadu s postroennymi obp2ektami pri fiksirovannyh znacheniyah $h$
budem rassmatrivatp1 ih pryamye (diskretnye ili nepreryvnye) summy po
vsem dopustimym $h$ (v chastosti,
$\Omega^{\SI}_{\BP;c}$ -- prostranstvo polubeskonechnyh form
Be1nksa--Peskina).

Otmetim, chto prostranstva differentsialp1nyh form Be1nksa--Peskina,
strunnyh polee0 Zigelya, kalibrovochno-invariantnyh strunnyh polee0,
va\-ku\-u\-mov Bovika--Radzhiva yavlyayut\-sya superprostranstvami, i
obp2ekty na nih takzhe yavlyayut\-sya superobp2ektami, no my dlya
kratkosti opuskaem pri\-stav\-ku ``super''. Formuly dlya dee0stviya
algebry Virasoro vo vseh uka\-zan\-nyh prostranstvah kak v ploskom,
tak i v iskrivlennom be1kgraunde, a takzhe dlya BRST-operatorov
soderzhat\-sya v [3,4]. Privedennye0 spisok dannyh polnostp1yu
harakterizuet svobodnuyu strunnuyu teoriyu polya, sa\-mo\-dee0\-st\-vie
trebuet vvedeniya dopolnitelp1nyh algebraicheskih struktur. Osnovnye
prostranstva teorii (prostranstva differentsialp1nyh form
Be1nksa--Peskina i strunnyh polee0 Zigelya) ne zavisyat ot metriki
be1k\-gra\-un\-da, kotoraya opredelyaet kak geometrodinamicheskie
obp2ekty strunnoe0 teorii polya: (psevdo)e1rmitovy metriki na
ukazannyh prostranstvah i BRST-operatory, tak i kalibrovochnye
harakteristiki -- strunnuyu svyaz\-nostp1 Gaussa--Manina i
sootvet\-stvuyushchie0 kovariantnye0 differentsial v
fokovsko-gostovom rassloenii. Otmetim, odnako, chto metriki i
BRST-operatory v prostranstvah differentsialp1nyh form
Be1nksa--Peskina i strunnyh polee0 Zigelya, rassmatrivaemyh kak
abstraktnye linee0nye pro\-st\-ran\-st\-va, nerazlichimy dlya razlichnyh
be1kgraundov (nezavisimostp1 teo\-rii vtorichnokvantovannyh
svobodnyh strun ot be1kgraunda), metrika na be1kgraunde
vosstanavlivaet\-sya pri rassmotrenii ukazannyh prostranstv kak
prostranstv geometricheskih obp2ektov na prostranstve vneshnih i
vnutrennih stepenee0 svobody struny. E1to pozvolyaet datp1
traditsionnuyu dlya strunnoe0 teorii polya interpretatsiyu metrik
na be1kgraunde i polee0 Yanga--Millsa kak ``nizkotemperaturnyh
predelov'' polee0 zamknutyh i ot\-kry\-tyh strun (hotya vozmozhny
alp1ternativy, naprimer, rassmatrivatp1 komponenty higgsovskih
polee0 dlya polee0 Yanga-Millsa Standartnoe0 Mo\-de\-li v tom zhe
sektore, chto i gravitatsionnoe, t.e. v sektore zamknutyh strun,
v silu sushchestvovaniya dopolnitelp1nyh nemetricheskih stepenee0
svo\-bo\-dy kogomologie0 algebry Virasoro s koe1ffitsientami v strunnyh
polyah [4,3]). Usloviya ploskosti svyaznosti Gaussa--Manina (ili,
chto to zhe sa\-moe, nilp1potentnosti kovariantnogo differentsiala)
sutp1 {\cyi strunnye urav\-ne\-niya E1e0nshtee0na\/} (perehodyashchie v
obychnye v ``nizkotemperaturnom pre\-de\-le''). Kak sledstvie, strunnye
uravneniya E1e0nshtee0na mozhno opredelitp1 i kak usloviya
sushchestvovaniya i nilp1potentnosti BRST-operatora Kato--Ogavy v
prostranstve kalibrovochno-invariantnyh strunnyh polee0 Zi\-ge\-lya.

Esli be1kgraund ne udovletvoryaet strunnym uravneniyam E1e0nshtee0na, to
vakuumy Bovika--Radzhiva ne sushchestvuyut, v e1toe0 situatsii estp1
ryad re\-tsep\-tov postroeniya strunnoe0 teorii polya, naprimer,
ispolp1zovanie ins\-tan\-to\-nov Bovika--Radzhiva [3]. Voprosy zavisimosti
(nezavisimosti) po\-lu\-cha\-yu\-shchih\-sya teorie0 ot be1kgraunda kak na
svyaznyh komponentah prostranstva reshenie0 strunnyh uravnenie0
E1e0nshtee0na, tak i v tselom, obsuzhdayut\-sya vo vtoroe0 iz rabot [3]
(i ssylki v nee0).

Dlya formulirovki teorii samodee0stvuyushchego strunnogo polya udobno
is\-polp1\-zo\-vatp1 opirayushchie0sya na idei nekommutativnoe0 geometrii
for\-ma\-lizm lokalp1nyh konformnyh polevyh algebr, podrobno izlozhennye0
v ori\-gi\-nalp1\-nyh rabotah [13] i obzore [14] (sm.takzhe [15]).

Rassmotirim, sleduya [3], prostranstvo $\Omega^{\cdot}_{\BP;\enl}=
\Omega^{\cdot}(\tilde\Bbb C^*,\Omega^{\SI}_{\BP;c})$ ras\-shi\-ren\-nyh
differentsialp1nyh form Be1nksa--Peskina (zdesp1 $\tilde\Bbb C^*$ --
uni\-ver\-salp1\-noe nakrytie prokolotoe0 v nule kompleksnoe0 ploskosti).
Opredelim tak\-zhe prostranstvo rasshirennyh strunnyh polee0 Zigelya
$\Omega^{\cdot}_{\sf;\enl}\!=\!
\Omega^{\cdot}(\tilde\Bbb C^*,(\Omega^{\SI}_{\BP;c})^*)$.
Formuly dlya dee0stviya algebry Virasoro v prostranstvah
$\Omega^{\cdot}_{\BP;\enl}$ i $\Omega^{\cdot}_{\sf;\enl}$ soderzhat\-sya
v [3]. Postroim takzhe rasshirennye BRST-operatory $Q_{\enl}$ i $Q^*_{\enl}$
kak vneshnie differentsialy v prostranstvah $\Omega^{\cdot}_{\BP;\enl}$ i
$\Omega^{\cdot}_{\sf;\enl}$ is\-ho\-dya iz BRST-operatorov $Q$ i $Q^*$.

\proclaim{\cyb Teorema 1 [3]} \cyi Prostranstvo $\Omega^{\cdot}_{\sf;\enl}$
dopuskaet strukturu BRST-in\-va\-ri\-ant\-noe0 lokalp1noe0 konformnoe0
polevoe0 algebry, kovariantnoe0 po ot\-no\-she\-niyu k strunnoe0 svyaznosti
Gaussa--Manina $\nabla_{\GM}$.
\endproclaim

\cyr
Takim obrazom, prostranstvo $\Omega^{\cdot}_{\sf;\enl}$ mozhet
rassmatrivatp1sya kak ne\-kom\-mu\-ta\-tiv\-nye0 kompleks de Rama (sr.[16,17])
po otnosheniyu k rasshirennym BRST-operatoram. E1tot kompleks
nazyvaet\-sya {\cyi rasshirennoe0 strunnoe0 po\-le\-voe0 algebroe0}.
Svyazp1 rasshirennoe0 strunnoe0 polevoe0 algebry s {\cyi
ne\-as\-so\-tsi\-a\-tiv\-noe0 strunnoe0 polevoe0 algebroe0
Arefp1evoe0--Volovicha}, realizuyushchee0sya v prostranstve strunnyh polee0
Zigelya i predstavlyayushchee0 soboe0 ne\-ko\-to\-ruyu snimayushchuyu
rasshirenie reduktsiyu assotsiativnoe0 rasshirennoe0 strun\-noe0 polevoe0
algebry, opisana v [3].

E1lementy rasshirennoe0 strunnoe0 polevoe0 algebry
$\Omega^{\cdot}_{\sf;\enl}$ obrazuyut al\-geb\-ru Li otnositelp1no
kommutatora. E1ta algebra Li dopuskaet tsent\-ralp1\-noe rasshirenie
pri pomoshchi mnimoe0 chasti (psevdo)e1rmitovoe0 met\-ri\-ki
$((\cdot|\cdot))$. Rassmotrim formy svyaznosti na $\tilde\Bbb C^*$ so
znacheniyami v $(\Omega^{\SI}_{\BP;c})^*$, t.e. kalibrovochnye polya na
$\tilde\Bbb C^*$ so znacheniyami v strunnyh polyah Zigelya; e1le\-men\-ty
$\Omega^0_{\sf;\enl}$ realizuyut infinitezimalp1nye ka\-lib\-ro\-voch\-nye
pre\-ob\-ra\-zo\-va\-niya e1tih polee0. E1ti kalibrovochnye preobrazovaniya
zamknuty (i tem samym my okazyvaemsya v situatsii strunnoe0 teorii
polya Vittena 1986 goda [5]), sootvet\-stvuyushchaya algebra Li
nazyvaet\-sya {\cyi vittenovskoe0 strun\-noe0 algebroe0 Li\/} i
oboznachaet\-sya $\cwit$ (krugovaya strelka $\circlearrowleft$
yavlyaet\-sya simvolom struny). Prostranstvo $\nabla_{\GM}$--kovariantnyh
e1lementov algebry Li $\cwit$ oboznachaet\-sya $\cwit_{\nabla^{\GM}}$ i
takzhe nazyvaet\-sya vittenovskoe0 strun\-noe0 algebroe0 Li.
Vittenovskaya strunnaya algebra Li $\cwit$ estp1 v tochnosti
tsen\-t\-ralp1\-noe rasshirenie kommutatornoe0 algebry nulevoe0 komponenty
ras\-shi\-ren\-noe0 strunnoe0 polevoe0 algebry, opisannoe vyshe.

Na prostranstve $\cwit^*$ (ili $\cwit^*_{\GM}$), dvoe0stvennom k
vittenovskoe0 strun\-noe0 algebre Li $\cwit$ (ili $\cwit_{\GM}$) zadany
kanonicheskie skobki Puassona -- skobki Li--Berezina, kotorye
kvantuyut\-sya kak takovye.

Skobki Li--Berezina v koprisoedinennom predstavlenii vittenovskoe0
strunnoe0 algebry Li mogut bytp1 redutsirovany do nepolinomialp1nyh
skobok v prostranstve funktsionalov na differentsialp1nyh formah
Be1n\-k\-sa--Peskina (ili vakuumah Bovika--Radzhiva), protsedura
gamilp1tonovoe0 reduktsii opisana v [3] i sleduet obshchee0 s\-heme
reduktsii skobok Li--Be\-re\-zi\-na (sm.naprimer, [18]). E1ti
nepolinomialp1nye skobki porozhdayut kva\-zi\-(psev\-do)\-al\-geb\-ru Li
(kvazialgebru v terminologii [19], sm.takzhe [20], i psevdoalgebru
v terminologii [18]) infinitezimalp1nyh ne\-po\-li\-no\-mi\-alp1\-nyh
kalibrovochnyh preobrazovanie0. E1ti nepolinomialp1nye
pre\-ob\-ra\-zo\-va\-niya na prostranstve vakuumov Bovika--Radzhiva byli
rassmotreny v [7]; oni obrazuyut kvazi(psevdo)algebru Li,
oboznachaemuyu $\czwie_{\nabla^{\GM}}$ i na\-zy\-va\-e\-muyu tsvibahovskoe0
strunnoe0 kvazi(psevdo)algebroe0 Li; so\-ot\-vet\-st\-vu\-yu\-shchaya
kvazi(psevdo)algebra Li na prostranstve differentsialp1nyh form
Be1nksa--Peskina oboznachaet\-sya $\czwie$ i imeet to zhe nazvanie.
Ne\-po\-li\-no\-mi\-alp1\-nye skobki Puassona realizuyut\-sya v funktsionalah
na pro\-st\-ran\-st\-ve $\czwie^*$ (ili $\czwie^*_{\nabla^{\GM}}$),
dvoe0stvennom k kvazi(psevdo)algebre Li $\czwie$ (ili
$\czwie_{\nabla^{\GM}}$). Tsvibahovskaya strunnaya
kvazi(psevdo)algebra Li mozhet bytp1 poluchena iz
neassotsiativnoe0 strunnoe0 polevoe0 algebry
Arefp1evoe0--Volovicha kak ee ``kommutatornaya'' algebra. Bolee
tochno, vys\-shie operatsii v mulp1tialgebre Sabinina--Miheeva [21],
postroennoe0 po tsvibahovskoe0 strunnoe0
kvazi(psevdo)algebre Li, sutp1 v tochnosti vys\-shie kommutatory
v neassotsiativnoe0 strunnoe0 polevoe0 algebre Arefp1evoe0--Volovicha.

Takim obrazom, nepolinomialp1naya strunnaya teoriya polya [7] v
pro\-st\-ran\-st\-ve differentsialp1nyh form Be1nksa--Peskina (ili vakuumov
Bovika--Radzhiva) mozhet bytp1 poluchena iz kubichnoe0 strunnoe0
teorii vit\-te\-nov\-s\-ko\-go tipa [5] v rasshirennom prostranstve pri
pomoshchi gamilp1tonovoe0 re\-duk\-tsii. Pri e1tom podhod rabot [7]
po nepolinomialp1noe0 teorii polya okazyvaet\-sya
e1kvivalentnym podhodu I.Ya.Arefp1evoe0 i I.V.Volovicha [6],
osnovyvayushchemusya na neassotsiativnoe0 strunnoe0 polevoe0 algebre.

\head\cyb 2. Kvantovogruppovaya struktura samodee0stvuyushchego strunnogo
polya\endhead

Zametim, chto kvantovanie samodee0stvuyushchego strunnogo polya
mozhet osushchestvlyatp1sya dvoyakim obrazom. Vo-pervyh, mozhno
kvantovatp1 sami ne\-po\-li\-no\-mi\-al1\-nye skobki Puassona, naprimer, v
formalizme asim\-p\-to\-ti\-ches\-ko\-go kvantovaniya [18]. Vo-vtoryh,
kvantovaya teoriya mozhet bytp1 poluchena pri pomoshchi kvantovoe0
reduktsii kvantovannyh skobok Li--Berezina na prostranstve
$\cwit_{\nabla^{\GM}}$; pri e1tom algebra kvantovyh nablyudaemyh
otozh\-des\-tv\-lya\-et\-sya s nekotoroe0 kvantovoe0 reduktsiee0
universalp1noe0 ober\-ty\-va\-yu\-shchee0 algebry
$\Cal U(\cwit_{\nabla^{\GM}})$
vittenovskoe0 strunnoe0 algebry Li. Oba va\-ri\-an\-ta predstavlyayut
nesomnennue0 interes s matematicheskoe0 tochki zreniya i
vyglyadyat dostatochno estestvennymi. Odnako, nizhe my postaraemsya
vy\-yavitp1 odin nemalovazhnye0 nyuans, znachitelp1no izmenyayushchie0
per\-vo\-na\-chalp1\-nye0 ``naivnye0'' vzglyad na predmet.

Otmetim, chto obp2ekty, postroennye vyshe, opisyvayut
sa\-mo\-dee0\-st\-vu\-yu\-shchuyu strunnuyu teoriyu polya tolp1ko na
drevesnom urovne, t.e. zadayut t.n. ``klassicheskuyu strunnuyu teoriyu
polya'' v terminologii [22]. Dlya togo chto\-by opisatp1 strunnuyu teoriyu
polya polnostp1yu i posledovatelp1no ne\-per\-tur\-ba\-tiv\-no, neobhodimo
ispolp1zovatp1 sleduyushchie0 klyuchevoe0 rezulp1tat.

\proclaim{\cyb Teorema 2 [3]} \cyi Vittenovskaya strunnaya algebra Li
$\cwit_{\nabla^{\GM}}$) (ili $\cwit$) dopuskaet strukturu bialgebry Li.
\endproclaim

\cyr
Dokazatelp1stvo e1toe0 teoremy opiraet\-sya na promezhutochnoe
ut\-ver\-zh\-de\-nie o tom, chto rasshirennaya strunnaya polevaya algebra
yavlyaet\-sya krossing-algebroe0 [3].

Itak, na kvaziklassicheskom urovne (sr.[23,18]) vyyavlena
kvan\-to\-vo\-grup\-po\-vaya struktura samodee0stvuyushchee0 strunnoe0
teorii polya. Takim obrazom, na kvantovom urovne [24] algebra
nablyudaenyh opisyvaet\-sya kvantovoe0 universalp1noe0
obertyvayushchee0 algebroe0 $\Cal U_q(\cwit_{\nabla^{\GM}})$ (ili
$\Cal U_q(\cwit)$), ili, tochnee, ee nekotoroe0 kvantovoe0
reduktsiee0, odnako, yavnaya konstruktsiya e1toe0 beskonechnomernoe0
algebry Hopfa neizvestna. Prezhde, chem peree0ti k opisaniyu
kvantovoe0 algebry nablyudaemyh samodee0stvuyushchee0 strunnoe0
teorii polya v konkretnyh sluchayah, razberem protsess reduktsii
na kva\-zi\-klas\-si\-ches\-kom urovne, inymi slovami, vyyasnim, vo chto
preobrazuet\-sya vittenovskaya strunnaya bialgebra Li pri reduktsii
skobok Li--Berezina, otvechayushchih vittenovskoe0 strunnoe0 algebre
Li, do nepolinomialp1nyh skobok Puassona (napomnim, chto sama
vittenovskaya strunnaya algebra Li perehodit v tsvibahovskuyu
strunnuyu kvazi(psevdo)algrbu Li).

\proclaim{\cyb Teorema 3}\cyi Tsvibahovskaya strunnaya kvazi(psevdo)algebra
Li $\czwie_{\nabla^{\GM}}$ (ili $\czwie$) obladaet strukturoe0
koyakobievoe0 kvazibialgebry.
\endproclaim

\cyr
Koyakobievy kvazibialgebry [25] yavlyayut\-sya klassom protobialgebr
Li, dvoe0stvennym yakobievym kvazibialgebram (kvazibialgebram Li v
ter\-mi\-no\-lo\-gii V.G.Drinfelp1da [26]).

Dlya dokazatelp1stva teoremy dostatochno primenitp1 reduktsiyu k
dublyu vittenovskoe0 strunnoe0 bialgebry Li s
translyatsionno--invariantnoe0 sko\-b\-koe0.

Napomnim, chto koyakobievy kvazibialgebry yavlyayut\-sya
in\-fi\-ni\-te\-zi\-malp1\-nym obp2ektom dlya puassonovyh kvazigrupp [25]
(v to vremya kak yakobievy kvazibialgebry -- dlya kvazipuassonovyh
grupp Li [26]). Takim obrazom, tsvibahovskaya strunnaya
kvazibialgebra Li realizuet na in\-fi\-ni\-te\-zi\-malp1\-nom urovne
kvaziklassicheskuyu versiyu nelinee0noe0 geometricheskoe0
al\-geb\-ry [27], chisto kvantovaya versiya kotoror0 do sih por
ne izvestna. Svyazp1 mezhdu strukturoe0 koyakobievoe0
kvazibialgebry i takim obp2ektom nelinee0noe0 geometricheskoe0
algebry kak mulp1tialgebra Sabinina--Mi\-he\-e\-va vyyavlena v [25].
Spetsifika beskonechnomernoe0 situatsii proyavlyaet\-sya v tom, chto,
po--vidimomu, ukazannym infinitezimalp1nym obp2ektam ne
sootvet\-stvuet nikakaya globalp1naya kvazigruppa.

Peree0dem teperp1 ot kvaziklassiki k yavnomu postroeniyu
(rasshirennoe0) kvantovoe0 algebry nablyudaemyh
$\Cal U_q(\cwit_{\nabla^{\GM}})$ v chastnyh sluchayah (i dlya
prostoty bez ucheta gostov, t.e. v fokovskom sektore).
Rassmotrim plos\-kie0 kompaktnye0 be1kgraund, izomorfnye0 faktoru
evklidova prostranstva po reshetke kornee0 poluprostoe0 algebry
Li $\frak g$. V e1tom sluchae ras\-shi\-ren\-naya strunnaya polevaya
algebra v prostranstve rasshirennyh ka\-lib\-ro\-voch\-no-in\-va\-ri\-ant\-nyh
strunnyh polee0 Zigelya $\Omega^0_{\sf;\enl}$ (tochnee, v ee
fokovskom sektore $\Omega^0_{\sf;\enl;F}$) yavlyaet\-sya
lokalp1noe0 konformnoe0 polevoe0 algebroe0, po\-lu\-chen\-noe0
perenormirovkoe0 potochechnogo proizvedeniya operatornyh
polee0 [3,14,15] (sm.takzhe [28]) iz algebry vershinnyh
operatorov, postroennyh po e1toe0 reshetke [29]. Sledovatelp1no,
linee0noe prostranstvo $\Omega^0_{\sf;\enl;F}$ mozhet bytp1
otozhdestvleno s prostranstvom $\Cal U(\hat\frak g_+))[[t]]$
polubeskonechnyh formalp1nyh stepennyh ryadov s koe1ffitsientami
v universalp1noe0 ober\-ty\-va\-yu\-shchee0 algebre polozhitelp1noe0
komponenty algebry Katsa--Mudi $\frak g_+$. Otmetim, chto
algebra vershinnyh operatorov porozhdaet\-sya svoimi tokami
(pervichnymi polyami spina 1), komponenty kotoryh obrazuyut
algebru Katsa--Mudi $\hat\frak g$, poe1tomu rasshirennaya
strunnaya polevaya algebra yavlyaet\-sya faktoralgebroe0
universalp1noe0 obertyvayushchee0 algebry $\Cal U(\hat\frak g)$
algebry Li $\hat\frak g$ po nekotoromu idealy $\Cal J$. Kak
sledstvie, vittenovskaya strunnaya algebra Li $\cwit_{\nabla^{\GM};F}$
(simvol $F$ oznachaet fokovskie0 sektor) yavlyaet\-sya
faktoralgebroe0 kommutatornoe0 algebry $\Cal
U_{[\cdot,\cdot]}(\hat\frak g)$ po idealu $\Cal J_{[\cdot,\cdot]}$.
Kvan\-to\-vaya versiya vittenovskoe0 strunnoe0 algebry Li
poluchaet\-sya sleduyushchim obrazom: rassmotrim kvantovuyu
universalp1nuyu obertyvayushchuyu $\Cal U_q(\hat\frak g)$,
snabzhennuyu $q$--kommutatorom; v silu sushchestvovaniya
$q$--verteksnoe0 kon\-s\-t\-ruk\-tsii dlya e1toe0 algebry, ideal
$\Cal J$ mozhet bytp1 prodeformirovan do ideala $\Cal J_q$
algebry $\Cal U_q(\hat\frak g)$, zamknutogo otnositelp1no
$q$--kommutatora; so\-ot\-no\-she\-niya mezhdu e1lementami
$\Cal U_q(\hat\frak g)/\Cal J_q$, zadavaemye $q$--kommutatorom,
i yavlyayut\-sya opredelyayushchimi v kvantovoe0 vittenovskoe0
strunnoe0 algebre $\Cal U_q(\cwit_{\nabla^{\GM};F})$.

Otmetim, chto v dannom klasse primerov, vazhnom s tochki zreniya
strun\-noe0 teorii, realizuyushchee0 kvantovuyu versiyu kalibrovochnyh
simmetrie0 strunnyh polee0, voznikal obp2ekt (rasshirennaya strunnaya
polevaya al\-geb\-ra), kotorye0 (kak i ego kvantovanie) estestvenno bylo
by rassmatrivatp1 v ramkah 2-petlevogo formalizma v teorii strun i
integriruemyh sistem [30]. Krome togo, bylo by interesno vyyavitp1,
kakim obrazom proyavlyaet\-sya modulyarnaya invariantnostp1
pervichnokvantovannoe0 struny v kvan\-to\-vo\-grup\-po\-vom formalizme
neperturbativnoe0 strunnoe0 teorii polya.

Itak, v dannoe0 rabote issledovany razlichnye kvantovogruppovye
struk\-tu\-ry strunnoe0 teorii polya kak na kvaziklassicheskom, tak i
na kvantovom urovne. Sformulirovany obshchie utverzhdeniya,
razobrany primery, pro\-de\-mon\-st\-ri\-ro\-va\-ny novye aspekty teorii
kvantovyh grupp (naprimer, svya\-zi s nelinee0noe0 geometricheskoe0
algebroe0 -- kvaziklassicheskie versii t.n. ``kvantovyh kvazigrupp''
i ``kvantovyh lup'') v kontekste vtorichnogo kvantovaniya struny
pri perehode s drevesnogo urovnya k po\-sle\-do\-va\-telp1\-no\-mu uchetu
petlevyh e1ffektov v ramkah neperturbativnoe0 strunnoe0 teorii polya.

\Refs\nofrills{\cyb Spisok literatury}
\roster
\item"[1]" {\cyie Grin M., Shvarts Dzh., Vitten E1.}, {\cyre Teoriya
superstrun. M., Mir, 1990.}
\item"[2]" {\cyie Morozov A.Yu., Perelomov A.M.} / {\cyre Sovrem.probl.matem.
Fundam.napravleniya. M., VINITI, 1989};
{\cyie Morozov A.Yu.} /\!/ {\cyre E1ChAYa. 1992. T.23(1). S.174-238; UFN. 1992.
T.162.}
\item"[3]" {\it Juriev D.} /\!/ Alg.Groups Geom. 1994. V.11. P.145-179
[e-version: hep-th/9403068]; Russian J.Math.Phys. 1996. V.4. P. 287-314;
J.Geom.Phys. 1995. V.16. P.275-300.
\item"[4]" {\it Juriev D.} /\!/ Lett.Math.Phys. 1991. V.22. P.1-6, 11-14;
1990. V.19. P.355-356; 1990. V.19. P.59-64.
\item"[5]" {\it Witten E.} /\!/ Nucl.Phys.B. 1986. V.268. P.253-294.
\item"[6]" {\it Aref'eva I.Ya., Volovich I.V.} /\!/ Phys.Lett.B. 1986.
V.182. P.159-163, 312-316, 1987. V.189. P.488.
\item"[7]" {\it Saadi M., Zwiebach B.} /\!/ Ann.Phys.(NY) 1989. V.192.
P.213-227; {\it Kugo T., Kunimoto H., Suehiro K.} /\!/ Phys.Lett.B. 1989.
V.226. P.48-54; {\it Sonoda H., Zwiebach B.} /\!/ Nucl.Phys.B. 1990. V.331.
P.592-628; {\it Kugo T., Suehiro K.} /\!/ Nucl.Phys.B. 1990. V.337. P.434-466;
{\it Zwiebach B.} /\!/ Commun.Math.Phys. 1991. V.136. P.83-118.
\item"[8]" {\it Volovich I.V.} /\!/ Class.Quant.Grav. 1987. V.4. P.83;
{\it Vladimirov V.S.} /\!/ Lett.Math.Phys. 1993. V.27. P.123.
\item"[9]" {\cyie Markov M.A.}, {\cyre Giperony i K-mezony. M., 1958.}
\item"[10]" {\it Segal G.} /\!/ Commun.Math.Phys. 1981. V.80.
P.301-342; {\it Kirillov A.A.} /\!/ Lect.Notes Math. 1984. V.970. P.49-67;
{\cyie Kirillov A.A., Yurp1ev D.V.} /\!/ {\cyre Funkts.anal.i ego pri\-lozh.
1986. T.20(4). S.79-80; 1987. T.21(4). S.35-46}; {\it Witten E.} /\!/
Commun.Math. Phys. 1988. V.114. P.1-53; {\it Kirillov A.A.} /
Infinite-dimensional Lie algebras and quantum field theory. World Scientific,
Teaneck, NJ, 1988, P.73-77; {\it Kirillov A.A., Juriev D.V.} /\!/
J.Geom.Phys. 1988. V.5. P.351-364; {\it Kirillov A.A.} / Operator algebras,
unitary rep\-re\-sen\-tations, envelopping algebras and invariant theory.
Birkhauser, Boston, 1992, P.73-83; {\it Juriev D.} /\!/ Adv.Soviet Math.
1991. V.2. P.233-247; {\it Kirillov A.A.} /\!/ Contemp.Math. 1993. V.145.
P.33-63.
\item"[11]" {\cyie Kirillov A.A.} /\!/ {\cyre Funkts.anal.i ego prilozh.
1987. T.21(2). S.42-45}.
\item"[12]" {\cyie Yurp1ev D.V.} /\!/ {\cyre Algebra i anal. 1990. T.2(2).
S.209-226}; {\it Juriev D.} /\!/ Russian J.Math.Phys. 1994. V.2(1). P.111-121.
\item"[13]" {\it Juriev D.} /\!/ Commun.Math.Phys. 1991. V.138. P.569-581,
1992. V.146. P.427; J.Funct. Anal. 1991. V.101. P.1-9.
\item"[14]" {\cyie Yurp1ev D.V.} /\!/ {\cyre UMN. 1991. T.46(4). S.115-138.}
\item"[15]" {\cyie Yurp1ev D.V.} /\!/ {\cyre TMF. 1994. T.101(3). S.331-348.}
\item"[16]" {\it Connes A.}, Introduction \`a la g\'eom\'etrie non commutative.
InterEditions, Paris, 1990.
\item"[17]" {\it Manin Yu.I.}, Topics in non-commutative geometry. Princeton
Univ.Press, Princeton, NJ, 1991.
\item"[18]" {\cyie Karasev M.V., Maslov V.P.}, {\cyre Nelinee0nye skobki
Puassona. Geo\-met\-riya i kvan\-to\-va\-nie. M., Nauka, 1991.}
\item"[19]" {\it Batalin I.} /\!/ J.Math.Phys. 1981. V.22. P.1837-1850.
\item"[20]" {\cyie Miheev P.O.} / {\cyre Nekotorye prilozheniya
differentsialp1noe0 geometrii. M., 1985, S.85-93 [VINITI: 4531-85Dep.]}; {\it
Mikheev P.O.} /\!/ Trans.Inst.Phys. Estonian Acad. Sci. 1990. V.66. P.54-66.
\item"[21]" {\cyie Sabinin L.V., Miheev P.O.} /\!/ {\cyre DAN SSSR. 1988.
T.297. S.801-805.}
\item"[22]" {\it Zwiebach B.} /\!/ Nucl.Phys.B. 1993. V.390. P.33-152.
\item"[23]" {\cyie Drinfelp1d V.G.} /\!/ {\cyre DAN SSSR. 1983. T.268(2).
S.285-287}; {\cyie Semenov-Tyan-Shanskie0 M.A.} /\!/ {\cyre Funkts.anal.i
ego prilozh. 1983. T.17(4). S.17-33}.
\item"[24]" {\cyie Drinfelp1d V.G.} /\!/ {\cyre Zap.nauchn.semin.LOMI. 1986.
T.155. S.19-49}; {\cyie Reshetihin N.Yu., Tahtadzhyan L.A., Faddeev L.D.}
/\!/ {\cyre Algebra i anal. 1989. T.1(2). S.178-206;} {\cyie Isaev A.P.} /\!/
{\cyre E1ChAYa. 1995. T.26(5). S.1204-1263}; {\cyie Zhelobenko D.P.},
{\cyre Pred\-stav\-le\-niya reduktivnyh algebr Li. M., Nauka, 1994.}
\item"[25]" {\it Bangoura M.} /\!/ C.R.Acad.Sci.Paris I. 1994. V.319.
P.975-978; {\it Bangoura M.\/}: Th\`ese l'Univ.Sci.Tech.Lille, no.1387, 1995.
\item"[26]" {\cyie Drinfelp1d V.G.} /\!/ {\cyre Algebra i anal. 1989. T.1(6).
S.114-148}; {\it Kosmann-Schwarzbach Y.} /\!/ C.R.Acad.Sci.Paris I. 1991.
V.312. P.391-394; {\it Kosmann-Schwarzbach Y.} / Ma\-the\-ma\-ti\-cal aspects
of field theory. Amer.Math.Soc., Providence, RI, 1992, pp.459-489.
\item"[27]" {\cyie Sabinin L.V.} / {\cyre Tkani i kvazigruppy. Kalinin
[Tverp1], 1988, S.32-37; Metody neassotsiativnoe0 algebry v
differentsialp1noe0 geometrii. Dobavlenie k russk. perev. Kobayasi S.,
Nomidzu K., Osnovy differentsialp1noe0 geometrii. T.1. M., Nauka,
1982}; {\cyie Sabinin L.V.} /\!/  {\cyre Trudy In-ta Matem. SO AN SSSR. 1989.
T.14. S.208-221}; {\cyie Miheev P.O., Sabinin L.V.} / {\cyre Probl.geometrii.
T.20. M., VINITI, 1988, S.75-100}; {\it Mikheev P.O., Sabinin L.V.} /
Quasigroups and loops. Theory and applications. Berlin, Heldermann Verlag,
1990, P.357-430;
{\it Sabinin L.V.} /\!/ Trans.Inst. Phys. Estonian Acad.Sci. 1990. V.66.
P.24-53.
\item"[28]" {\cyie Yurp1ev D.V.} /\!/ {\cyre TMF. 1992. T.93(1). S.32-38.}
\item"[29]" {\it Frenkel I., Lepowsky J., Meurman A.}, Vertex operator
algebras and the Monster. Acad. Press, New York, 1988; {\it Frenkel I.B.,
Huang Y.-Z., Lepowsky J.} /\!/ Memoires AMS. 1993. V.494; {\it Dong C.,
Lepowsky J.}, Generalized vertex algebras and relative vertex operators.
Birkh\"auser, 1993; {\it Kac V.G.}, Vertex operator algebras for beginners.
Birkh\"auser, 1996.
\item"[30]" {\cyie Morozov A.Yu.} /\!/ {\cyre UFN. 1994. T.164(1). S.1.}
\endroster
\endRefs
\enddocument